\documentclass[10pt]{article}

\usepackage{amsmath,amscd,amsfonts,amsthm, graphicx, xcolor}

\newcommand{\shrinkmargins}[1]{
  \addtolength{\textheight}{#1\topmargin}
  \addtolength{\textheight}{#1\topmargin}
  \addtolength{\textwidth}{#1\oddsidemargin}
  \addtolength{\textwidth}{#1\evensidemargin}
  \addtolength{\topmargin}{-#1\topmargin}
  \addtolength{\oddsidemargin}{-#1\oddsidemargin}
  \addtolength{\evensidemargin}{-#1\evensidemargin}
  }

\shrinkmargins{.7}

\newcommand{\field}[1]{\mathbb{#1}}

\newcommand{\Z}{\field{Z}}

\newcommand{\R}{\field{R}}

\newcommand{\ra}{\rightarrow}

\newcommand{\set}[1]{\{#1\}}

\newcommand{\beq}{\begin{displaymath}}
\newcommand{\eeq}{\end{displaymath}}
\newcommand{\beqn}{\begin{equation}}
\newcommand{\eeqn}{\end{equation}}
\newcommand{\del}{\delta}
\newcommand{\eps}{\epsilon}
\newcommand{\calT}{\mathcal{T}}

\theoremstyle{plain}
\newtheorem{thm}{Theorem}
\newtheorem{prop}[thm]{Proposition}

\theoremstyle{definition}

\newtheorem{exmp}[thm]{Example}

\newtheorem{prob}[thm]{Problem}

\theoremstyle{remark}

\title{Convergence rates for ordinal embedding}
\author{Jordan S. Ellenberg and Lalit Jain}
\date{April 2019}
\begin{document}

\maketitle

\begin{abstract} We prove optimal bounds for the convergence rate of ordinal embedding (also known as non-metric multidimensional scaling) in the $1$-dimensional case.  The examples witnessing optimality of our bounds arise from a result in additive number theory on sets of integers with no three-term arithmetic progressions.  We also carry out some computational experiments aimed at developing a sense of what the convergence rate for ordinal embedding might look like in higher dimensions.
\end{abstract}

\section{Introduction}
\label{s:intro}

The problem of {\em non-metric multidimensional scaling}, also known as {\em ordinal embedding}, is the following:  given an integer $n$,  a metric space $M$ (often a low-dimensional manifold) with distance function $d$, and a set $\Sigma$ of ordered quadruples $(i,j,k,l) \in [1 \ldots n]^4$ we aim to find an embedding of $[1 \ldots n]$ into $M$, which we think of as an $n$-tuple of points $(x_1,x_2, \ldots, x_n)$,  satisfying the constraints
\beq
d(x_i,x_j) < d(x_k,x_l)
\eeq
for all $(x_i, x_j, x_k, x_l) \in \Sigma$. 

This problem is of special relevance to the analysis of data obtained from human respondents.  One can try to measure dissimilarity between concepts quantitatively; for instance, one might ask subjects ``On a scale of $1$ to $7$, how similar are the colors red and purple?"  But this creates an obvious calibration problem:  one person's $2$ might be another person's $4$.  The question ``which pair is more similar: red and purple or green and black?" is likely to yield more consistent answers. 

In the present paper, we consider the special case where every quadruple in $\Sigma$ is of the form $(i,j,i,k)$.  That is, each element of $\Sigma$ specifies an answer to the question ``which of the two points $x_j$ or $x_k$ is closer to $x_i$?" 
%
%That is, we impose either
%\beq
%d(x_i,x_j) < d(x_i,x_k)
%\eeq
%or
%\beq
%d(x_i,x_k) < d(x_i,x_j).
%\eeq

Every constraint of this form is called a ``triplet comparison''.  We concentrate on this case because these triplet comparisons seem to be psychologically natural and reliable when gathering data from human participants; it is easier to say whether red is more like purple or more like green than it is to say whether red is more like green or a cow is more like a fish \cite{borggroenen}.

%We now face a natural question:  to what extent can we determine the $N$-tuple $x_1, \ldots, x_N \in M$ from a set of triplet constraints $\Sigma$?

There is a critical difference between non-metric multidimensional scaling and classical MDS, in which distances between points are specified exactly rather than by comparison.  Take $M=\R^d$.  If all pairwise distances betweeen $n$ points in $\R^d$ are known, it is easy to see that the locations of the points themselves are specified, up to a distance-preserving affine linear transformation of $\R^d$.  However, even if we have {\em all} possible non-metric data -- that is, if we know for all $(i,j,k)$ whether $x_i$ is closer to $x_j$ or to $x_k$ --  it is clear that the location of the points is not specified.  If $(x_1, \ldots, x_n)$ satisfies all the constraints in $\Sigma$, then typically so do at least some sufficiently small perturbations of $(x_1, \ldots, x_n)$.

The question of convergence for ordinal embedding thus splits naturally into two parts.
\begin{itemize}
\item How many triplet comparisons do we need in order to determine the answer to all $n {n -1 \choose 2}$ triplet comparisons?
\item If all triplet comparisons are specified, within what error can we determine $(x_1, \ldots, x_n)$?
\end{itemize}

In this paper we restrict our attention to the second question.  (For more on the first question, which can be seen as a non-metric version of {\em rigidity}, see \cite{robkevin} and \cite{jainjamiesonnowak}.)   Following the notation of \cite{ac}, we say that a function $f: M \ra N$, where $M$ and $N$ are metric spaces, is {\em weakly isotonic} if, for every three points $m,m',m'' \in M$, we have that
\beq
 d_M(m,m') < d_M(m,m'')
\eeq
if and only if
\beq
d_N(f(m),f(m')) < d_N(f(m),f(m'')).
\eeq  

That is, $M$ and $f(M) \subset N$ agree on all triplet comparisons.  By a slight abuse of notation, we say that two $n$-tuples $x=(x_1, \ldots, x_n)$ and $y = (y_1, \ldots, y_n)$ in an ambient metric space $M$ are weakly isotonic if the function from the induced metric spaces $\set{x_1, \ldots, x_n}$ to $\set{y_1, \ldots, y_n}$ sending $x_i$ to $y_i$ is weakly isotonic. 

%We use $x$ to refer either to the ordered $n$-tuple of points $(x_1, \ldots, x_n)$ or to the discrete subspace $\set{x_1, \ldots, x_n} \subset M$; context will make it clear which is meant.

We begin by showing that weakly isotonic subsets of $[0,1]$ cannot differ by too much, as long as at least one of the sets has no long gaps.  In order to state this result precisely, we need to introduce some notation.  If $x$ is a finite subset of $M$, the {\em Hausdorff distance} between $x$ and $M$, denoted $\delta_H(x,M)$, is the smallest $\alpha$ such that every point of $M$ is within $\alpha$ of a point of $x$.  In particular, a finite subset $S$ of $[0,1]$ is close in Hausdorff distance to $[0,1]$ if there is no large interval disjoint from $S$.

If $x,y \in M^n$ are two $n$-tuples we denote by $d_\infty(x,y)$ the maximum value of $d_M(x_i,y_i)$ over $i \in [1,\ldots,n]$.  By a {\em similarity} of a metric space $M$ we mean a homeomorphism from $M$ to $M$ which multiplies all distances by a scalar.  A similarity of Euclidean space is just the composition of a dilation with an isometry.  Of course, if $x$ and $y$ are related to each other by a similarity, then $x$ and $y$ are isotonic, despite possibly being far away from each other in $d_\infty$.  It is thus natural to consider $\min_A d_\infty(x,Ay)$ as a measure of the distance between $x$ and $y$.  (Though we do not primarily adopt this point of view in the present paper, one might also consider an $L_p$ norm of the $d_M(x_i,y_i)$ instead of the maximum, or $L_\infty$ norm; we denote the resulting distance by $d_p(x,y)$.)

We are now ready to state our first result, which concerns the distance between weakly isotonic pairs of subsets of the interval.

\begin{prop} Suppose $x=(x_1, \ldots, x_n)$ and $y = (y_1, \ldots, y_n)$ are subsets of $[0,1]$ such that $\del_H(x,[0,1]) \leq \alpha$.  Suppose $x$ and $y$ are weakly isotonic.  Then there exists a similarity $A$ of $\R$ such that  $d_\infty(x,Ay) = O_\eps(\alpha^{1-\eps})$\footnote{The notation $O_\epsilon(\cdot)$ indicates that the big-O is hiding a constant that may depend on $\epsilon$.} for all $\epsilon > 0$.
\label{pr:intervalbound}
\end{prop}

A result of Kleindessner and von Luxburg~\cite{KvL} shows, under the same conditions, that $d_\infty(x,y)$ goes to $0$ as $\alpha$ does. Though they do not state convergence rates explicitly, their proof in the case of $M=[0,1]$ gives $d_\infty(x,y) = O(\alpha^{1/2})$.  

%Note that the condition on the diameter of $y$ is acting to avoid the trivial case where $y$ is a dilation of $x$; in that case, $y$ and $x$ agree on all triplet comparisons but may be far apart.

The condition on the Hausdorff distance between $x$ and $[0,1]$ may seem a bit artificial at first, but some version of this condition is necessary, as the following example demonstrates.

\begin{exmp} Let $x_1, \ldots, x_{n/2}$ lie in $[-0.1,0.1]$, and let $x_{n/2+1}, \ldots, x_n$ lie in $[0.9,1.1]$.  Define $y_i$ to be $x_i$ for $i \in 1,\ldots,N/2$ and $x_i + 1000$ for $i \in N/2+1, \ldots, N$.  It is easy to see that $x$ and $y$, despite being very far apart,  agree on all triplet comparisons.  In fact, it is clear that $x$ and $y$ can be arbitrarily far apart while agreeing on all triplets.
\label{ex:clusters}
\end{exmp}

In some sense, the lesson of Example~\ref{ex:clusters} is that we can't expect the triplet comparisons to correctly locate a set of points in a metric space $M$ if we've chosen the wrong metric space $M$ in the first place.  In the example, the problem is that the points naturally cleave into two separate clusters, which suggests that we should have been trying to embed them into two disjoint copies of the interval, not into a single interval.  Suppose, for instance, that we were asking subjects to rate the similarity of words, and that the words were actually half colors and half animals.  The subjects will reliably and correctly report that red is more like blue than red is like bird; but an embedding of the words in Euclidean space commits us to a judgment that bird is, say, {\em ten times} as far from red as red is from blue.  It seems clear that we shouldn't expect or desire our method to specify measurements of this kind.

\bigskip

Our next result is that Proposition~\ref{pr:intervalbound} is optimal up to epsilon exponents:

\begin{prop}  For all sufficiently small $\alpha$, there exist subsets $x=(x_1, \ldots, x_n)$ and $y = (y_1, \ldots, y_n)$ of $[0,1]$ such that
\begin{itemize}
\item $\del_H(x,[0,1]) \leq \alpha$;
%\item  $y$ has diameter $1$; and
\item  $x$ and $y$ are weakly isotonic
\end{itemize}
with $d_\infty(x,Ay) = \Omega_\eps(\alpha^{1+\eps})$ for all similarities $A$.
\label{pr:ap}
\end{prop}

The novelty here is that the proof of Proposition~\ref{pr:intervalbound}, given in the next section, rests on number-theoretic results about large subsets of $\Z$ containing no three-term arithmetic progressions.  To extend this result to higher dimensions, one would need bounds on the number of dense subsets of $[0,1]^d$ containing no approximately isosceles triangles; we formalize this question below as Problem~\ref{prob:isosceles}.  This problem seems potentially of interest in its own right as a question in combinatorial geometry.

The situation in higher dimensions is less settled.  Arias-Castro (see the remark after \cite[Theorem 3]{ac}) proves a convergence rate for the situation where all quadruple comparisons, not merely triple comparisons are known, and says a similar result can be obtained for triple comparisons: if $U$ is a bounded, connected, open domain in $\R^d$ satisfying some modest conditions, and $x$ is an $n$-tuple at Hausdorff distance at most $\alpha$ from $U$, and $y$ is weakly isotonic to $x$, then $d_\infty(x,Ay) = O(\alpha^{1/2})$ for some similarity $A$ of $\R^d$.  As we have shown in Proposition~\ref{pr:intervalbound}, this is not optimal even for $M = [0,1]$.  For large $d$, the situation is murkier still.  It's not hard to show that, if $x$ is a set of $n$ points chosen uniformly in $[0,1]^d$, the Hausdorff distance between $x$ and $[0,1]^d$ is at most $n^{-1/d+\eps}$, and that at least one point of $x$ can be moved a distance $\sim n^{-(2+\epsilon)}$ without changing any triplet comparisons.  The idea is as follows: for each triple $x_i,x_j,x_k$ let $d(i,j,k)$ be the least difference between any two edgelengths of the triangle formed by $x_i,x_j,x_k.$  One can check that each of $x_i,x_j,x_k$ can be moved a distance $d(i,j,k)$ in any direction without changing $x_i,x_j,x_k$.  The probability that any individual $d(i,j,k)$ is less than $n^{-2.01}$ is on order $n^{-2.01}$, so the total number of triples $x_i,x_j,x_k$ with $d(i,j,k) < n^{-2.01}$ has expected value $\sim n^{0.99}$, which means that with high probability there exists an $x_i$ (indeed many $x_i$) which is not contained in any such triple, and which can be thus moved a distance $n^{-2.01}$ without changing any triplet comparisons.

Thus there exist $x,y$ which are weakly isotonic, with $\del_H(x,[0,1]^d) \sim \alpha$ and $d_\infty(x,y)$ at least on order $\alpha^{2d+\epsilon}$.

%We provide a lower bound for the optimal error rate in in the case of $[0,1]^d$ with its $L^1$ metric.
%
%\begin{prop}  Let $M$ be $[0,1]^d$ endowed with the $L^1$ metric.  There exist subsets $x=(x_1, \ldots, x_n)$ and $y = (y_1, \ldots, y_n)$ of $M$ such that $\del_H(x,M) \leq \alpha$, $y$ has diameter $1$, and $x$ and $y$ are weakly isotonic, with $d_\infty(x,y) = \Omega_\eps(\alpha^{2d+\eps}).$
%\label{pr:sidon}
%\end{prop}
%
%\begin{rem}  Our intuition is that there should be an example in $[0,1]^d$ with the Euclidean metric which yields a bound of the same order as Proposition~\ref{pr:sidon}, but we do not see an easy way to adapt the $L^1$ example to work in $L^2$.
%\end{rem}

What is the optimal convergence rate for subsets close in Hausdorff distance to $[0,1]^d$?  In particular, does the rate improve with $d$, like the lower bound above, or is it independent of $d$, like Arias-Castro's upper bound in \cite{ac}?  Some mild evidence in favor of the former guess is provided by the following proposition, which shows that a {\em random} perturbation of $x$ is highly unlikely to be weakly isotonic to $x$, unless the perturbation is very small.

\begin{prop}  Let $x = (x_1, \ldots, x_n)$ be a subset of $[0,1]^d$ of size $n$.  Let $y = (y_1, \ldots, y_n)$ be a subset of $\R^d$ in which $y_i$ is chosen uniformly at random from the Euclidean ball of size $\beta>n^{-1}$ around $x_i$.  Then the probability that $y$ is weakly isotonic to $x$ is bounded above by $\exp(-cn)$ for some constant $c>0$. 
\label{pr:random}
\end{prop}

In particular, the Hausdorff distance $\alpha = \del_H(x,[0,1]^d)$ between this random subset and the cube is likely to be on order $n^{-1/d}$; thus we can say that a random perturbation of $x$ of size $\alpha^d$ is very unlikely to be weakly isotonic to $x$.

It is conceptually helpful to think of these results geometrically.  If $x = (x_1, \ldots, x_n)$ is a subset of $[0,1]^d$, we can write $P_x$ for the subset of $([0,1]^d)^n$ consisting of all $y$ which are weakly isotonic to $x$.  We know almost nothing about $P_x$, essentially because the triplet conditions cutting out $P_x$ are given by non-definite quadratic forms and hence are non-convex.  Thus we don't know that $P_x$ is convex, or even that it is connected.  Propositions~\ref{pr:intervalbound} and \ref{pr:ap} should be thought of as saying something about the {\em radius} of $P_x$, i.e. the maximal distance of any point in $P_x$ from $x$.\footnote{More precise would be to say the maximal distance from the orbit of $x$ under similarities.}  If $P_x$ tends to be ``round'', we might expect its radius to behave similarly to the size of the largest ball around $x$ it contains.  If $P_x$ tends to be ``pointy'', one might expect the radius to be substantially larger. Proposition~\ref{pr:random} says that, for every $x$, a ball around $x$ of radius larger than $n^{-1}$ lies mostly outside $P_x$.  The argument presented in the paragraph following Proposition~\ref{pr:ap} says that, for a randomly chosen $x$, the region $P_x$ is unlikely to be entirely contained within a ball of radius $n^{-2} = \alpha^{2d}$ around $x$ (i.e. the radius is likely to be at least $n^{-2}$ - though the proposition does not imply that $P_x$ contains a ball of this radius either).  We note that in Proposition~\ref{pr:random} the controlling parameter is the number of points, not the Hausdorff distance from $[0,1]^d$.  It is instructive in this context to consider $x$ as in Example~\ref{ex:clusters}; for such an $x$, the region $P_x$ has radius of {\em constant} size, but is still very long and thin, with most of the ball of radius $n^{-1}$ around $x$ outside $P_x$ as guaranteed by Proposition~\ref{pr:random}.

In the two sections that follow, we prove all three propositions we've stated.  We conclude in the final section by presenting some experimental data which gives hints as to the truth about the precision to which triplet comparisons specify an embedding.

\subsection*{Acknowledgments} The first author is supported by NSF grant DMS-1700884 and by a Simons Foundation Fellowship.  The authors are grateful for the hospitality of the Institute for Foundations of Data Science at UW-Madison where some of this work was carried out. Thanks to Laura Balzano for her feedback on a draft.

\section{Upper bounds:  proof of Proposition~\ref{pr:intervalbound}}

We will prove a slightly more precise version of Proposition~\ref{pr:intervalbound}.  Our proof is very much in the spirit of Kleindessner and von Luxburg's \cite[Lemma 9]{KvL}.

\begin{prop} Let $x = (x_1, \ldots, x_n)$ be an subset of the interval $[0,1]$ containing both $0$ and $1$, whose Hausdorff distance from $[0,1]$ is at most $\alpha$.  Let $y = (y_1, \ldots, y_n)$ be a subset of $[0,1]$ which is weakly isotonic to $x$.  Then there is a similarity $A$ of $\R$ such that $|x_i - A y_i| < 2 \alpha (\log_2 (1/\alpha) + 3/2)$ for all $i$. % In particular, if $A$ is the similarity $z \ra \lambda z + b$, we have  $d_\infty(x,Ay) = O(\alpha^{1 - \epsilon})$.
\label{th:interval}
\end{prop}

We remark that the condition that $x$ contains $0$ and $1$ is not meaningfully restrictive, since if $x$ had smaller diameter we could make that condition hold by applying a similarity to $x$.  This could slightly increase $d_H(x,[0,1])$ from $\alpha$ to a quantity of the form $\alpha + o(\alpha)$; we have imposed the condition $0,1 \in x$ in order to avoid writing this down.

\begin{proof}

%First of all, we can and do apply a similarity to $x$ until it has diameter $1$; in fact, we may further assume that the $x_i$ are in increasing order with $x_1 = 0$ and $x_n = 1$.   

We note that in $[0,1]$ the notion of ``betweenness" is determined by triplet comparisons:  $x_k$ is between $x_i$ and $x_j$ precisely when $d(x_i,x_j) \geq d(x_i,x_k)$ and $d(x_i,x_j) \geq d(x_j,x_k)$.  It follows that $y_1, \ldots, y_n$ are either in the same order as $x_1, \ldots, x_n$ or the opposite order.  Reflecting $y$ around $1/2$ if necessary, we may assume they are in the same order, and applying a dilation, we can assume that $y_1= 0$ and $y_n=1$.

Let $x_{1/2}$ be the smallest element of $x$ which is greater than $1/2$. By the condition on Hausdorff distance, $x_{1/2}$ lies in $[1/2,1/2+2\alpha]$.  Write $y_{1/2}$ for the point of $y$ corresponding to $x_{1/2}$ (that is, it occurs in the same ordinal position among the $y_i$ as $x_{1/2}$ does among the $x_i$.)  Since $x_{1/2}$ is closer to $1$ than to $0$, so is $y_{1/2}$; that is, $y_{1/2} \geq 1/2$, so in particular, 
\beq
x_{1/2} - y_{1/2} \leq 2\alpha.
\eeq
Now let $x_{1/4}$ be the smallest element of $x$ greater than $(1/2)(0 + x_{1/2})$, and let $x_{3/4}$ the smallest element of $x$ greater than $(1/2)(x_{1/2} + 1)$.  Then $x_{1/4}$ must lie in the range $[1/4,1/4 + 3 \alpha]$, while $x_{3/4}$ lies in $[3/4, 3/4 + 3 \alpha]$.  If $y_{1/4}$ is the element of $y$ corresponding to $x_{1/4}$, then $y_{1/4}$ must exceed the mean of $0$ and $y_{1/2}$, which is at least $1/4$.  Similarly, $y_{3/4} > 3/4$.

Continue in this manner, defining $x_{m/2^k}$, for each odd $m$ in $[1,2^k-1]$, to be the smallest element of $x$ greater than the mean of $x_{(m-1)/2^k}$ and $x_{(m+1)/2^k}$.  We always have
\beq
x_{m/2^k} \in [m/2^k, m/2^k + a_k \alpha]
\eeq
where $a_k$ is the increasing sequence defined by the recurrence
\beq
a_0 = 0, a_1 = 2, a_k = (1/2)(a_{k-1} + a_{k-2}) + 2
\eeq
and is bounded above by $2k$ (in fact, $a_k$ is asymptotic to $(4/3)k$).  Moreover, we have
\beq
y_{m/2^k} > m/2^k.
\eeq

Now let $k$ be $\lceil \log_2 1/\alpha \rceil$.  Choose some $x_i$ in $x$.  The fact that $2^{-k} \leq \alpha$ implies that there is some integer $M$ such that $M/2^k$ lies in the range $[x_i - (a_k+1) \alpha, x_i - a_k \alpha]$.  Write $M/2^k$ in lowest terms as $m / 2^\ell$ with $m$ odd and $\ell \leq k$. 

%If $m$ is odd and $\ell \leq k$, then $x_{m/2^\ell} <= m/2^\ell + 2 \ell \alpha$.  

%Choose $m/2^\ell$ to be maximal among all fractions with denominator dividing $2^k$ such that
%\begin{equation}
%m/2^\ell + a_\ell \alpha < x_i.
%\label{eq:bound2ml}
%\end{equation}
%The fact that $2^{-k} \leq \alpha$ implies that there is some integer $m_0/2^\ell$ in the range $(x_i - a_k \alpha, x_i - (a_k+1)\alpha)$, and any such $m_0/2^\ell$ satisfies \eqref{eq:bound2ml}; so $m/2^\ell$ is at least  $x_i - (a_k+1)\alpha$.

Since
\beq
x_i > m/2^\ell + a_k \alpha \geq m/2^\ell + a_\ell \alpha \geq x_{m/2^\ell}
\eeq
we also have $y_i > y_{m/2^\ell} \geq m/2^\ell$.  Thus
\beq
x_i - y_i \leq (a_k+1)\alpha.
\eeq
for all $i$. We can run the whole argument upside down (that is, replacing $x_i$ by $1-x_i$ and $y_i$ by $1-y_i$) yielding
\beq
x_i - y_i \leq - (a_k+1)\alpha
\eeq
for all $i$.  Now $a_k \leq 2k$ and $k \leq \log_2 (1/\alpha) + 1$, which gives the desired result.

\end{proof}

\section{Lower bounds: proofs of Propositions~\ref{pr:ap} and \ref{pr:random}}

The convergence rate of $\alpha^{1-\epsilon}$ proved in Proposition~\ref{pr:intervalbound} is nearly optimal.  It turns out that constructing subsets of the interval which are difficult to approximate by triplet comparisons is very close to the well-known problem in additive number theory of constructing large subsets of $\set{1,2,\ldots,N}$ which contain no three integers in arithmetic progression.

In particular:  by a result of Graham~\cite{graham:vanderwaerden}, for all positive integers $k$ there exists a subset $S$ of $\Z$ such that 

\begin{itemize}
\item $S$ is contained in the interval $[1,M]$, with $M \geq k^{c \log k}$ for some absolute constant $c$;
\item $S$ has no three terms in arithmetic progression;
\item $S$ has no gaps between successive terms of size greater than $k$.
\end{itemize}

Now let $x = (x_1, \ldots, x_{|S|})$ be the set of points $\set{s/M: s \in S} \in [0,1]$.

We note that $X$ has no gaps of size greater than $k/M$; in other words, the Hausdorff distance between $x$ and $[0,1]$ is at most $k/2M$, which quantity we denote by $\alpha$.  Moreover, for every triplet $x_i, x_j, x_k$, we have that
\beq
2x_k - x_i - x_j \neq 0
\eeq
(because $S$ has no three terms in arithmetic progression) and thus
\beq
|2x_k - x_i - x_j| \geq 1/M
\eeq
(because every $x_i$ lies in $(1/M)\Z$.)

In particular, this means that if $y = (y_1, \ldots, y_{|S|})$ satisfies
\beq
|y_i  - x_i| < 1/2M
\eeq
for all $i$, then $y$ is weakly isotonic to $x$.  Now $k$ grows more slowly than any power of $M$, so $1/M  = \Omega(\alpha^{1+\epsilon})$ as desired.

We need to show not only that $d_\infty(x,y) = \Omega(\alpha^{1+\epsilon})$ but that the same is true of $d_\infty(x,Ay)$ for any similarity $A$.  We argue as follows.  Choose two consecutive elements $x_i$ and $x_{i+1}$ of $x$ with $|x_i - x_{i+1}|$ on order $\alpha$, and take $y_i = x_i - \beta$ and $y_{i+1} = x_{i+1} + \beta$, for some $\beta$ on order $\alpha^{1+\epsilon}$.  Then $d_\infty(x,y) \geq \beta$.  Now choose another consecutive pair $x_j,x_{j+1}$ satisfying the same conditions and again take $y_j = x_j-\beta, y_{j+1} = x_{j+1} + \beta$. For all $k\neq i,i+1,j,j+1$ let $y_{k} = x_k$.

Let $A$ be a similarity of $\R$. For $A$ to satisfy $d_\infty(x,Ay) < \beta$, we would need $Ay_i > y_i$ and $Ay_{i+1} < y_{i+1}$.  This in turn implies that the unique fixed point of $A$ lies between $x_i$ and $x_{i+1}$.  The same reasoning applies to $x_j,x_{j+1}$; since the fixed point of $A$ cannot lie inside both intervals, we conclude that $d_\infty(x,Ay) \geq \beta$.

We remark that we have shown more than we need to; in order for $d_\infty(x,y)$ to be large we only need {\em one} $x_i$ to have substantial freedom of motion, while in this case we have shown that {\em all} $x_i$ do.  So one might think of Proposition~\ref{pr:ap} as giving lower bounds for $d_p(x,y)$ for all $p$, not only $p=\infty.$

\medskip

We turn now to Proposition~\ref{pr:random}, whose proof is quite simple.  Choose some $x_i$ in $x$.  The other $n-1$ points of $x$ are all at bounded distance from $x_i$, so we can find $x_j,x_k$ which form a ``nearly isosceles triangle" with $x_i$; that is,
\beq
|d(x_i,x_j) - d(x_i, x_k)|  < C n^{-1}
\eeq
for some $C$.  This being the case, a random perturbation of $x$ of size on order $n^{-1}$ will modify $x_i,x_j,x_k$ by independent real numbers of size around $n^{-1}$, and thus changes the sign of $d(x_i,x_j) - d(x_i, x_k)$ with some positive probability $p$.  That is, with probability $p$ the perturbation changes one of the triplet comparisons.  Now choose a new $x_{i'}$ from the $n-3$ elements of $x$ other than $x_i, x_j, x_k$, and choose $x_{j'}, x_{k'}$ so that $x_{i'},x_{j'}, x_{k'}$ is nearly isosceles, as before.  The probability that our perturbation leaves the sign of $d(x_{i'},x_{j'}) - d(x_{i'}, x_{k'})$ unchanged is again $1-p$, and this event is independent of the change of sign of $d(x_i,x_j) - d(x_i, x_k)$.  Repeating this process $\gamma n$ times, for some small constant $\gamma > 0$, we find that the probability that the perturbation of $x$ leaves all triplet comparisons unchanged is at most $\exp(-cn)$, as claimed.

\section{Experimental results}

\begin{figure}
    \centering
    \includegraphics[width=\linewidth]{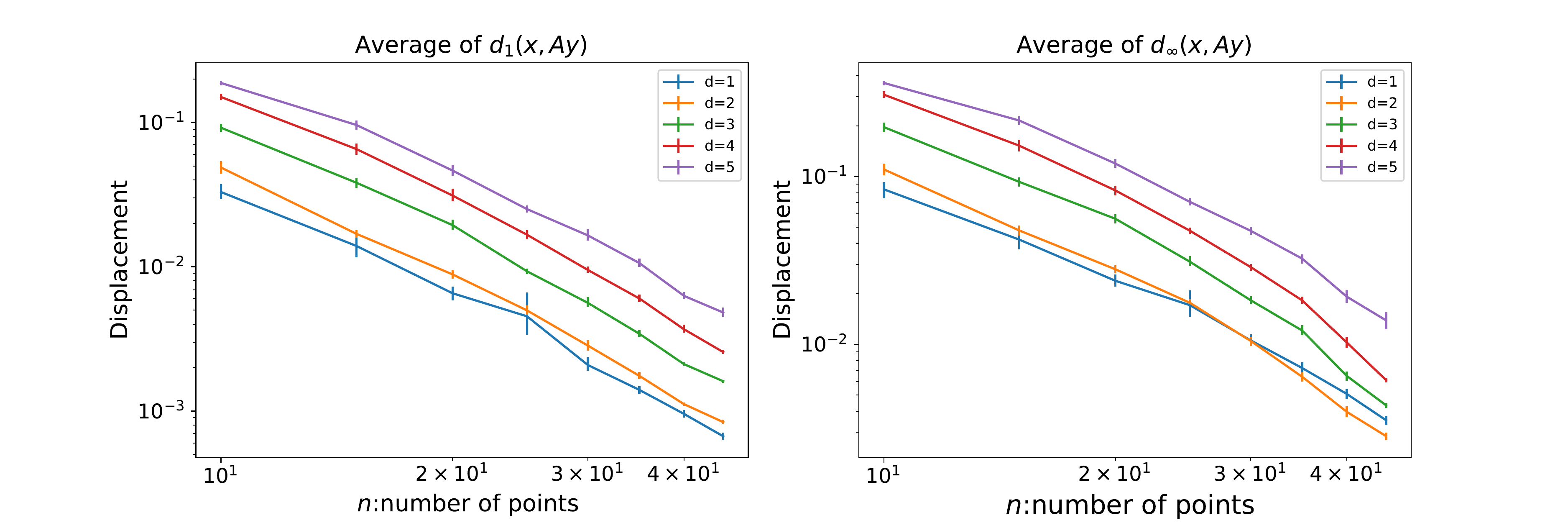}
    \caption{Average $\ell^{1}$ displacement and $\ell^{\infty}$ displacement.}
    \label{fig:maxavg}
\end{figure}

In this section, we describe some computational experiments investigating the convergence rates of non-metric multidimensional scaling in various dimensions.  Let $x = x_1, \ldots, x_n$ be points of $\R^d$.  We record all triplet comparisons arising from $x_1, \ldots, x_n$, denoted $\calT_x$ and then, starting from scratch, attempt to find another set of points $y = y_1, \ldots, y_n$ in $\R^d$ which agrees with $x$ on all triplet comparisons; that is, we try to solve the quadratic satisfaction problem of generating a $y$ which is weakly isotonic to $x$.  Having done so, we find the minimum (as described below an approximation) of $d_\infty(x,Ay)$ over all similarities $A$. For each dimension $d\in \{1,2,3,4,5\}$, choice of $n=\{10,15,20, 25, 30, 35, 40, 45\}$, we carried out this procedure for 100 different sets of points uniformly chosen from $[0,1]^{d}$ and averaged the results. In Figure~\ref{fig:maxavg} we show the main result for the main focus of this paper, the average maximum displacement,  $d_\infty(x,y) =  \sup_i d(x_i,y_i)$ and also the {\em average} displacement $d_1(x,y) = N^{-1} \sum_i d(x_i,y_i)$. We now discuss what this figures show and describe how our data was generated. 
%$\{x^{i}_{n,d}\}_{i=1}^{100}, x^{i}_{n,d}\subset [0,1]^d, |x^{i}_{n,d}| = n$

In general, the quadratic problem of finding a set of points $y$ satisfying the constraints in $\calT_{x}$ is a difficult problem and there are many attempted solutions, all involving nonconvex optimization techniques, that date back to the early days of psychometrics \cite{kruskalshepard, borggroenen}. Unfortunately, none of these methods have theoretical justification, but in practice will find a set of points satisfying all the constraints. For our empirical studies, we solved the following optimization problem:
\begin{equation} \label{eq:opt}
\min_{y} \sum_{ i,j,k\in \calT_x} \max(0, 1-(\|y_j - y_i\|_2^2 - \|y_k - y_i\|_2^2))
\end{equation}
The loss function inside the sum encourages $\|y_j - y_i\|_2^2 - \|y_k - y_i\|_2^2>0$ and is commonly used in constraint satisfaction problems \cite{gnmds, bowerjainbalzano}. The resulting optimization is easily solved to a local optima using stochastic gradient descent methods. 

\begin{figure}
    \centering
    \includegraphics[width=\linewidth]{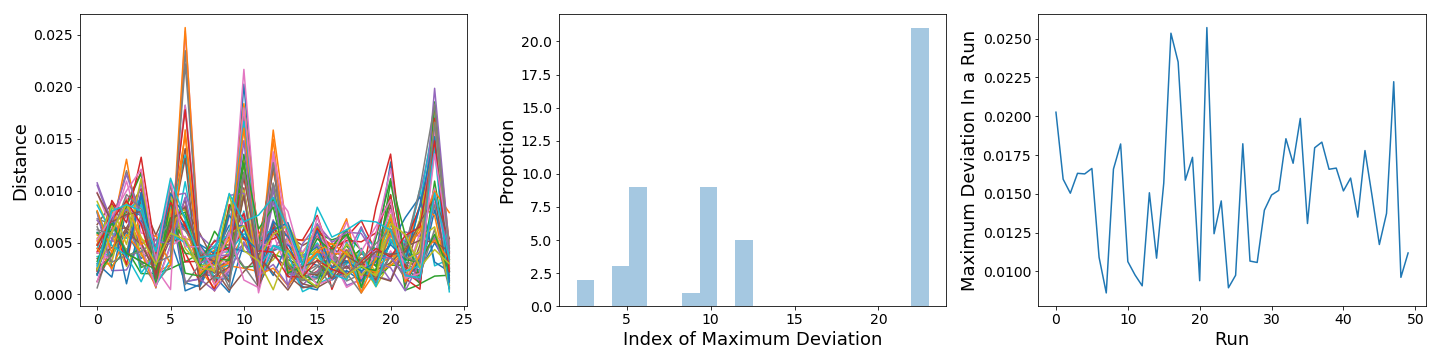}
    \caption{For each of 50 different embeddings of a set of points $x$, the displacement of each point (after the Procrustes similarity is applied) is shown. The second plot is the frequency of how many times each point attains the maximum displacement in any of these runs, and the final plot is the maximum displacement in each run.}
    \label{fig:runsrepeated}
\end{figure}

For each set of points $x$ chosen as described above, the algorithm in equation \ref{eq:opt} was used to construct a new embedding $y$ satisfying all the triplets in $\calT_x$. In general, $y\not\subset [0,1]^d$, so an orthogonal Procrustes transformation, i.e. the similarity $A$ minimizing $d_2(x, Ay) := \sum_i d(x_i, y_i)^2 $, was found using standard orthogonal Procrustes analysis \cite{gower75}. In general we can not expect the minimizer of $d_2(x, Ay)$ to give rise to a minimizer of $d_{\infty}(x, Ay)$. However, since 
\[\inf_{A} d_{\infty}(x, Ay) \leq \inf_{A} d_{2}(x, Ay) \]
this certainly provides an upper bound. Note that individual runs of the embedding algorithm/optimization in equation \ref{eq:opt} can result in very different embeddings. For example, Figure \ref{fig:runsrepeated} the displacement of each of 25 points in 2 dimensions over 50 different re-embeddings is shown. As can be seen, on any given run the index of the point which is displaced the most can vary, but it tends to be a single point (in the example, index 20). However the maximum displacement of a point on any given run was on index 6 - an extremal configuration we may not have found if we only computed one embedding. To mitigate this issue, for each $x$, 10 different embeddings $y$ were found by solving the optimization problem \ref{eq:opt}. For each $x$, out of the ten embeddings, the $y$ (after the Procrustes similarity was applied) with the point that had maximum deviation from its original value, was chosen. The resulting average of $d_1(x,Ay)$ and $d_{\infty}(x, Ay)$ was then computed over all 100 choices of $x$ (along with appropriate bootstrap confidence intervals). Figure \ref{fig:maxavg} shows the result.

%The theoretical results of this paper concern the maximum displacement.  Experimentally, we also considered the     %Again, for each of the pairs $(x,y)$ computed above, we find the minimum of $d_1(x,Ay)$ over all similarities $A$.

If we choose $n$ points $x_1, \ldots x_n$ uniformly at random in $[0,1]^d$, the Hausdorff distance between $x$ and $[0,1]^d$ is typically on order $n^{-1/d}$.  Thus, \cite{ac} would guarantee that $\min_A d_\infty(x,Ay)$ is at most $n^{-1/2d}$.  On the other hand, Proposition~\ref{pr:random} and the discussion proceeding it suggest that $d_\infty(x,Ay)$ might typically behave more like $n^{-c}$ for some $c$ independent of $d$. 

% COMMENT: WE are really saying something serious about the method we are using. Our thing is not consistent with the cell having small diameter. Snapshot of the cell around a point. 
The curves plotted in Figure \ref{fig:maxavg} appear to have slope in the range $[-2,-2.5]$ in the log-log plot, and in particular this slope does not appear to be changing with $d$.  In particular, the decay appears to be faster than the rate of $n^{-1/2d}$ provided by \cite{ac}. From the discussion before Proposition 3 we know that with high probability the maximum displacement of a set of $n$ uniformly chosen points in $[0,1]^d$ is {\em at least} $n^{-2}$. We take the results recorded in Figure~\ref{fig:maxavg} to be very modest evidence that $n^{-2}$ is the correct rate of decay.

One reason for caution:  as in the description at the end of Section~\ref{s:intro}, write $P_x$ for the region in $([0,1]^d)^n$ parametrizing $y$ weakly isotonic to $x$.  We are interested in the radius of $P_x$; what our algorithm is doing is (repeatedly) finding a point $y$ in $P_x$, typically on the boundary of $P_x$, and computing the distance between $y$ and $x$.  If $P_x$ is typically ''round," then this is likely a good estimate for the radius of $P_x$.  If $P_x$ is typically ''pointy," then the question becomes:  is the algorithm likely to converge to the ``pointy part" of $P_x$, which may be rather far from $x$?  If not, it is quite possible that the algorithm here systematically underestimates the radius of $P_x$ by finding $y \in P_x$ whose distance from $x$ is highly submaximal.

It is possible that more extensive experimentation along the lines of that recorded in Figure~\ref{fig:runsrepeated} would be useful for teasing this out.  Repeated runs of the algorithm with the same $x$ amounts to sampling from the $P_x$, albeit with no provable control on the distribution.  Suppose $P_x$ were very pointy, or even disconnected with one component containing $x$ and one far away. One might imagine that, under these circumstances, some of the lines in Figure~\ref{fig:runsrepeated} would be well above the bulk of the chart, corresponding to those runs yielding a $y$ in the faraway portion of $P_x$.  We have not observed any behavior like this in practice, but as yet we have no way to rule it out.  What's more, since we are using $x$ uniformly chosen in $[0,1]^d$, it may be that $P_x$ is round for some $x$ but highly pointy for a small set of exceptional $x$; indeed, we {\em know} this is the case if we don't place a constraint on Hausdorff distance, as Example~\ref{ex:clusters} shows.

\begin{figure}
    \centering
    \includegraphics[width=\linewidth]{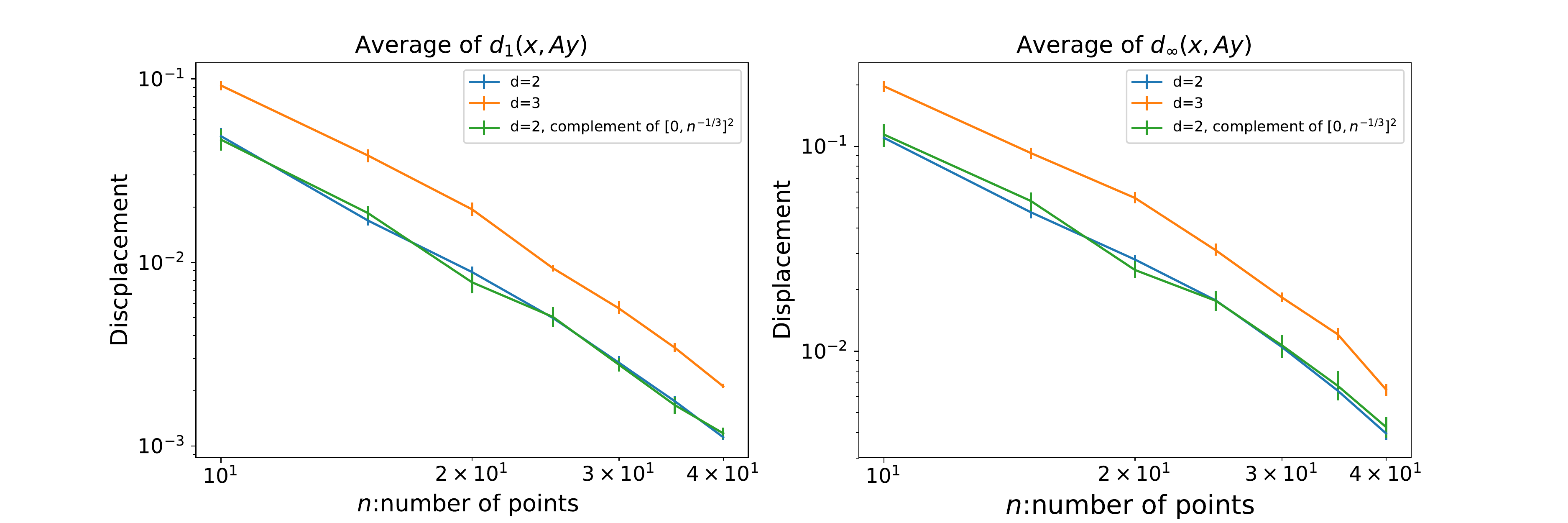}
    \caption{Average $\ell^1$ and $\ell^\infty$ displacements for sets of points in $[0,1]^2$, $[0,1]^3$ and outside of the box $[0,n^{-1/3}]$}
    \label{fig:box}
\end{figure}

With this issue in mind, we report on one more experiment:  we draw $n$ points from the complement of a square of side length $n^{-1/3}$ in $[0,1]^d$.  We now have a subset of $[0,1]^d$ of size $n$ whose Hausdorff distance $\alpha$ from the square is not $n^{-1/2}$, as would be the case for uniformly chosen points, but $n^{-1/3}$.  So one can ask: does the error in reconstruction of $x$ behave roughly like $n^{-2}$, or more like $\alpha^{2d} = n^{-4/3}$? Figure \ref{fig:box} seems the imply the former is much closer to the truth, and perhaps evidence that for uniformly chosen points perhaps a typical $P_x$ is indeed roughly round with radius $n^{-2}$ independent of the Hausdorff dimension.

%However note that this could be due to the fact that we have only looked at a few small values of $n$, and a more pronounced effect could be observed for larger $n$.  

Finally:  we emphasize that in Proposition~\ref{pr:ap} we construct, in the case $d=1$, examples of $x \subset [0,1]$ for which $P_x$ has radius on order $n^{-1-\epsilon}$, and even contains a ball of that radius.  Can we construct similar $x$ in $[0,1]^d$?  Of course we could simply place $n$ points on a line in $[0,1]^d$ and achieve the same convergence rate; to avoid degenerate examples like this, we ask that the set of points converges in Hausdorff distance to $[0,1]^d$.  The sets in Proposition~\ref{pr:ap} are obtained by means of subsets of the line with very few approximate 3-term arithmetic progressions.  A 3-term arithmetic progression is nothing more than an isosceles triangle contained in a line, so in higher dimensions the goal is to avoid approximate isosceles triangles.  For $\beta>0$ a real parameter, we say a triangle is {\em $\beta$-isosceles} if some vertex is within $\beta$ of the line equidistant between the other two vertices.

\begin{prob}
Let $S$ be a subset of $[0,1]^d$ containing no $\beta$-isosceles triangle.  What is the minimum Hausdorff distance between $S$ and $[0,1]^d$?
\label{prob:isosceles}
\end{prob}

Good bounds for Problem~\ref{prob:isosceles} would provide good lower bounds for convergence rate of ordinal embedding in higher dimensions. In the spirit of Proposition~\ref{pr:ap}, one might, for instance, ask: is there a subset of the $N \times N$ integer grid of size $N^{2-\epsilon}$ whose complement contains no large squares and which has no three points forming an isosceles triangle? 

%We note that if $S$ has size $n$, the distances to any given element $s \in S$ from the remaining elements in $S$ are $\sim n$ numbers in a bounded interval, so there must be two that differ by at most $\sim n^{-1}.$  We thus see that $n \leq \beta^{-1}$, which means that $d_H(x,[0,1]^d) \geq \beta^{1/d}$.  Can this bound be achieved?  For example, 

%{\bf OLD STUFF, TO BE COMMENTED OUT}

%{\bf For Lalit:}  What I would like to see here is two experiments.  First, take $n$ random points in $[0,1]$, compute all triplet comparisons, and use your algorithm to find a set of $n$ points in $[0,1]$ satisfying these comparisons.  My guess is that the output of the algorithm would place each point within $n^{-1}$ or so of the ``correct" spot -- in fact I wouldn't be surprised if starting from random points you end up within $n^{-2}$ a lot of the time.  It also might depend on whether we are talking about "maximum amount any point moves," or "average amount each point moves."

%Then, do the same thing with $n^2$ points in $[0,1]^2$.  (So the Hausdorff distance between pointset and manifold should again be roughly $1/n$.)   Here my guess is that we would see output at worst $n^{-2}$ and maybe on average more like $n^{-4}$ from the truth.


\begin{thebibliography}{1}



\bibitem[AL07]{gnmds}{Agarwal, Sameer, Josh Wills, Lawrence Cayton, Gert Lanckriet, David Kriegman, and Serge Belongie. ``Generalized non-metric multidimensional scaling.'' In Artificial Intelligence and Statistics, pp. 11-18. 2007.}

\bibitem[A15]{ac} Arias-Castro, Ery, Some theory for ordinal embedding, arXiv preprint arXiv:1501.02861 (2015).

\bibitem[BG03]{borggroenen}{Borg, Ingwer, and Patrick Groenen. ``Modern multidimensional scaling: Theory and applications.'' Journal of Educational Measurement 40.3 (2003): 277-280.}

\bibitem[BJB17]{bowerjainbalzano}{Bower, Amanda, Lalit Jain, and Laura Balzano. ``The Landscape of Non-Convex Quadratic Feasibility.'' In 2018 IEEE International Conference on Acoustics, Speech and Signal Processing (ICASSP), pp. 3974-3978. IEEE, 2018.}

\bibitem[G06]{graham:vanderwaerden} Graham, Ron. "On the growth of a Van der Waerden-like function." {\em Integers} 6 (2006).


\bibitem[G75]{gower75}{Gower, John C. "Generalized procrustes analysis." Psychometrika 40, no. 1 (1975): 33-51.}

\bibitem[KvL14]{KvL} Kleindessner, Matth{\"a}us and von Luxburg, Ulrike, Uniqueness of ordinal embedding, {\em Proceedings of COLT}, 40--67 (2014).

\bibitem[JJN15]{jainjamiesonnowak}{Jain, Lalit, Kevin G. Jamieson, and Rob Nowak. "Finite sample prediction and recovery bounds for ordinal embedding." {\em Advances In Neural Information Processing Systems.} 2016.}

\bibitem[JN11]{robkevin} Jamieson, Kevin G. and Nowak, Robert D, Low-dimensional embedding using adaptively selected ordinal data, {\em Proceedings of Communication, Control, and Computing (Allerton)}, 1077--1085 (2011). 


\bibitem[KS74]{kruskalshepard}{Kruskal, Joseph B., and Roger N. Shepard. "A nonmetric variety of linear factor analysis." Psychometrika 39, no. 2 (1974): 123-157.}



\end{thebibliography}
\end{document}